\documentclass[12pt]{amsart} 
\usepackage[curve,v2]{xy}

\usepackage{amssymb,latexsym,amsfonts, amsthm, amsmath}

\newcommand{\tensor}{\otimes}

\newcommand{\Ext}{\operatorname{Ext}}
\newcommand{\codim}{\operatorname{codim}}
\newcommand{\pdim}{\operatorname{proj-dim}}
\newcommand{\reg}{\operatorname{reg}}

\newcommand{\h}{\operatorname{h}}

\newcommand{\wt}[1]{\widetilde{#1}}

\newcommand{\wtp}{\widetilde{\P}^n}

\newcommand{\wts}{\widetilde{\Sigma}}

\newcommand{\ses}[3]{0\rightarrow#1\rightarrow#2
   \rightarrow#3\rightarrow0}

\newcommand{\PP}{\ensuremath{\mathbb{P}}}
\newcommand{\ZZ}{\ensuremath{\mathbb{Z}}}

\newcommand{\mO}{\ensuremath{\mathcal{O}}}

\newcommand{\mI}{\ensuremath{\mathcal{I}}}

\newcommand{\F}{{\mathcal F}}

\newcommand{\I}{{\mathcal I}}

\renewcommand{\P}{{\mathbb{P}}}

\newcommand{\Z}{{\mathbb{Z}}}

\newcommand{\K}{{\mathbb{K}}}

\newcommand{\of}{\, \mbox{\small{$\circ$}} \:}

\newtheorem{thm}{Theorem}[section]   
\newtheorem{cor}[thm]{Corollary}     
\newtheorem{lemma}[thm]{Lemma}         
\newtheorem{prop}[thm]{Proposition}  

\theoremstyle{definition} 
\newtheorem{conj}[thm]{Conjecture}        
   
\newtheorem{ex}[thm]{Example}        
\newtheorem{rmk}[thm]{Remark}

\newenvironment{rem}[2]{\refstepcounter{thm} \label{#2} 
\par \medskip \noindent {\bf #1 \thethm .}}{\par \medskip}

\begin{document}

\pagenumbering{arabic}

\title{Syzygies of the secant variety of a curve}

\author[Jessica Sidman]{Jessica Sidman}

\address{415a Clapp Lab,
Department of Mathematics and Statistics,
Mount Holyoke College,
South Hadley MA 01075}

\email{jsidman@mtholyoke.edu}

\author[Peter Vermeire]{Peter Vermeire}

\address{Department of Mathematics, 214 Pearce, Central Michigan
University, Mount Pleasant MI 48859}

\email{verme1pj@cmich.edu}


\date{\today}

\begin{abstract} 
We show the secant variety of a linearly normal smooth curve of degree at least $2g+3$ is arithmetically Cohen-Macaulay, and we use this information to study the graded Betti numbers of the secant variety.
\end{abstract}

\maketitle

\section{Introduction}
We work throughout over an algebraically closed field $\K$ of characteristic $0$.  A well-known result dating back to Castelnuovo states that if $C \subset \P^n$ is a linearly normal curve of genus $g$ with $\deg C\geq 2g+1,$ then $C$ is projectively normal and hence is arithmetically Cohen-Macaulay (ACM).  Our main result is
\begin{thm}
\label{intromainACM}
If $C \subset \P^n$ is a smooth linearly normal curve of genus $g$ and degree $d \geq 2g+3$, then its secant variety $\Sigma$ is ACM. 
\end{thm}

Using the Auslander-Buschbaum theorem \cite[\S 19]{eisenbud}, this tells us that a minimal free resolution of the homogeneous coordinate ring of $\Sigma,$ has length equal to $\codim \Sigma$, and the remainder of this paper is devoted to studying the syzygies among the defining equations of $\Sigma$.  One can get a rough idea of the behavior of the syzygies of a coherent sheaf using Castelnuovo-Mumford regularity, which may be defined in terms of cohomology.   Recall that a coherent sheaf $\F$ on $\P^n$ is $m$-regular (in the sense of Castelnuovo and Mumford)  if $H^i(\P^n, \F(m-i)) = 0$ for all $i > 0,$ and that the regularity of $\F$ is the infimum over all $m$ so that it is $m$-regular.  In Corollary~\ref{h^4} we show that if $\reg \I_{\Sigma}<5,$ then $C$ is rational and $\reg \I_{\Sigma} = 3.$

To describe our results on syzygies more precisely, we set up some notation.  Let $S = \K[x_0, \ldots, x_n].$  Any finitely generated graded $S$-module $M$ has a minimal free resolution
\[
0 \to \oplus S(-j)^{\beta_{r,j}} \to \cdots \to  \oplus S(-j)^{\beta_{1,j}} \to  \oplus S(-j)^{\beta_{0,j}}  \to M \to 0,
\]
where the \emph{graded Betti numbers} $\beta_{i,j}$ are uniquely determined by minimality.  It is convenient to display the $\beta_{i,j}$ in a graded Betti diagram in which the $(i,j)$ entry is $\beta_{i, i+j}.$
\[
\begin{array}{c|cccc}
           & 0 & 1 & 2 & 3 \\
             \hline
       0& \beta_{0,0} & \beta_{1,1}& \cdots \\
          1& \beta_{0,1}& \beta_{1,2}& \cdots \\
          2& \beta_{0,2} & \beta_{1,3}& \cdots\\
\end{array}
\]  As in \cite{E2} we say that the Betti numbers $\beta_{i, i+k}$ in the $i$th row of the Betti diagram form the degree $k+1$ \emph{linear strand} if $M = S/I$ for some homogeneous ideal $I.$ In this case, $\beta_{1, k+1}$ is the number of minimal generators of $I$ in degree $k+1.$   We say that $M$ is $m$-regular if $\beta_{i,j}=0$ for all $j > i+m.$ 

If $C \subset \P^n$ is a linearly normal curve of genus $g$ and degree $d \geq 2g+3,$ we obtain several results as consequences of the Cohen-Macaulay condition.    We give explicit formulas for several graded Betti numbers in Corollary~\ref{duality} and  Proposition~\ref{cubics}, showing that 
\begin{itemize}
\item $\beta_{1,3} = \binom{n+1}{3}-(d-2)n-3g+1$
\item $\beta_{2,4} = \beta_{1,4}+\beta_{1,3}(n+1)-\binom{n+4}{n}+P_{\Sigma}(4)$ 
\item $\beta_{n-3, n+1}= \binom{g+1}{2},$
\end{itemize}
where $P_{\Sigma}(x)$ is the Hilbert polynomial of $\Sigma.$  Note that via Theorem~\ref{intromainACM} there are exactly $n-3$ syzygy modules in the resolution of $S_{\Sigma},$ and if $g\geq1,$ then Corollary~\ref{h^4} implies that the last module is generated by elements of degree $\leq n+1.$  Thus, $\beta_{n-3, n+1}$ is the bottom right corner of the graded Betti diagram, and it depends only on the genus of the curve.

We compute the Hilbert polynomial of $S_{\Sigma}=S/I_{\Sigma}$ by relating it to the Hilbert polynomial of a curve of degree $D$ and genus $G$ gotten by intersecting $\Sigma$ with a plane of codimension 2.  
\begin{thm} \label{thm: hilbert poly} The Hilbert polynomial of $S_{\Sigma}$  agrees with its Hilbert function for all positive integers and is given by
\[
D\binom{m+2}{3}+(1-G)\binom{m+1}{2}+\alpha_1m+\alpha_0,
\]
where $\alpha_1 = \binom{n+2}{2}-(n+1)-3D-2(1-G)$ and $\alpha_0 = -\binom{n+2}{2}+2(n+1)+2D+1-G.$  
\end{thm}

We also obtain a nonvanishing result on the graded Betti numbers of higher secant varieties.

\begin{thm}\label{beyond p}
Let $C$ be a smooth curve of genus $g$ embedded into $\P^n$ via a line bundle $L$ of degree $d \geq 2g+2k+p+1$ and $\Sigma_k$ be its variety of secant $k$-planes.  Suppose that $L = L_1 \otimes L_2$ where $|L_1|=s \leq |L_2|=t.$  If $s+1 \geq k+2,$ then the length of the degree $k+2$ linear strand of $S_{\Sigma_k}$ is at least $s+t-2k-1.$  In particular, if $L$ is a general line bundle of degree $d \geq 2g+2k+p+1,$ then $\beta_{s+t-2k-1, s+t-k}(\Sigma_k) \neq 0.$
\end{thm}

We briefly sketch part of the picture of what is known about syzygies of high degree curves to put our results in context.  The homogeneous coordinate ring of a curve of degree at least $2g+1$ is 1-regular if $g=0$ and has regularity two otherwise.  A variety satisfies property $N_0$ if it is projectively normal, satisfies $N_1$ if its ideal is generated by quadrics and satisfies $N_p$ for $p\geq 2$ if all syzygies are linear through the $p$th stage of the resolution.  Through work of Green and Lazarsfeld \cite{mgreen},\cite{GL1},\cite{GL2} we know that if $d \geq 2g+1+p,$ then the curve satisfies $N_p.$    Moreover, from Theorem 8.17 in \cite{E2} due to Schreyer we know that $\beta_{p+\lfloor \frac{g}2 \rfloor,p+\lfloor \frac{g}2 \rfloor+1} \neq 0.$  Furthermore, as a consequence of duality, the ``last'' graded Betti number is $\beta_{n-1, n+1} = g.$  (See Chapter 8 of \cite{E2} for a nice discussion.)

We now extend and refine the conjectures in \cite{vermeiresecreg}.  
\begin{conj}\label{mainconj}
Suppose that $C \subset \P^n$ is a smooth linearly normal curve of genus $g$ and degree $d \geq 2g+2k+1+p$, where $p,k\geq0.$ Then 
\begin{enumerate}
\item $\Sigma_k$ is ACM and has regularity $(2k+2)$ unless $g=0,$ in which case the regularity is $k+1$.  
\item $\beta_{n-2k-1, n+1} = \binom{g+k}{k+1}.$ 
\item $\Sigma_k$ satisfies $N_{k+2,p}$, where a variety $X$ satisfies $N_{d,p}$ if its ideal is generated by elements of degree $d$ and all syzygies are linear through the $p$th stage of the resolution, as defined in \cite{eghp}.
\end{enumerate}
\end{conj}
As described above, the full conjecture is known to hold for $k=0$.  Further, by \cite{vBH} and \cite{fisher} it holds for $g\leq1$.  In this work, we show that parts (1) and (2) hold for $k=1$.  After the completion of this work, it was shown in \cite{n3p} that part (3) also holds for $k=1$.

We illustrate the behavior that we have seen with the example below  
\begin{ex}\label{intro g=2}
At the suggestion of David Eisenbud we used ideas of Frank-Olaf Schreyer to compute the ideal of a genus 2 curve embedded in $\P^7$.  Let $\overline{C}$ be a plane curve of degree $5$ with $4$ nodes.  If we blow up the four nodes in $\P^2$ and consider the linear system $|5H-4\Sigma E_i|$, the restriction of this system to the proper transform $C$ of $\overline{C}$ has degree $9=2g+5$, and embeds $C\subset\P^7$ as a smooth curve of genus $2$. 
Using  \emph{Macaulay 2} \cite{M2}  we can write down a basis of $|5H-4\Sigma E_i|$   over the rationals.  We  display the Betti diagram of the coordinate ring, where ``-'' denotes a zero entry.  
\[
\begin{array}{c|ccccccc}
           & 0 & 1 & 2 & 3 & 4 & 5 & 6\\
             \hline
       0& 1 & - & - & - & - & - & -\\
          1& - & 19& 58& 75& 44& 5 & -\\
          2& - & - &-  & - & - & 6 & 2\\
\end{array}
\]
Note that the quadratic strand of the resolution has length 5 but that the curve satisfies $N_4$ but not $N_5.$

Using code developed for \cite{sidmanSullivant}, we computed the ideal of $\Sigma.$  From the Betti diagram we see that the cubic strand of the resolution has length 2 and that $\beta_{4,8} = 3$ as predicted by Conjecture~\ref{mainconj}.
\[
\begin{array}{c|ccccc}
          &  0 & 1 & 2 & 3 & 4\\
           \hline
          0& 1  & -  & - & - & -\\
          1& - & - & - & - & -\\
          2& - & 12& 16 & - & -\\
          3& -  & -  & - &4 & -\\
          4& -  & -  & - & 4 &3\\
\end{array}
\]
Comparing the diagram to the statement of Corollary~\ref{cor: tail}, we see that the three unknowns at the tail of the resolution are all zero here as in Example~\ref{g=2}.
\qed
\end{ex}
 
We give a brief outline of the structure of the paper.  The ACM condition is treated in \S 3.  To understand the ACM condition, we work geometrically to show that cohomology groups vanish.  The key observation is that there is a desingularization $\wts \to \Sigma$ such that $\wts$ is a $\P^1$-bundle over the symmteric square of $C,$ which we denote by $S^2C$, and hence the cohomology of the structure sheaf of $\wts$ is the same as that of $S^2C$, which is easier to understand.  As $\Sigma$ has non-rational singularities, the higher direct image sheaves of the ideal of $\wts$ do not vanish, but there is another divisor whose ideal sheaf has the same direct image and whose higher direct images do vanish. (See Lemma \ref{tech}.)    Making the exact relationships between these objects precise is the bulk of our work.  The technical preliminaries are summarized in \S 2.  
We examine the graded Betti diagram of $S_{\Sigma}$ in \S 4.  

To improve readability we have written out some arguments which are surely well-known to experts, but are perhaps not easily available in the standard references.
\subsection*{Acknowledgements}  We thank David Eisenbud for his suggestions, and we thank two anonymous referees for suggestions which improved this paper.  We would not have discovered the statements of our main theorems without help from \emph{Macaulay 2} \cite{M2}.  In computing secant equations we used code developed with the help of Mike Stillman in conjunction with the first author's work with Seth Sullivant on \cite{sidmanSullivant}.  The first author is partially supported by NSF grant DMS 0600471 and the Clare Boothe Luce Program.

\section{Setup and notation}
Suppose that $X \subset \P^n$ is a variety.  We let $\mO_X$ and $\I_X$ denote the structure sheaf and ideal sheaf of $X.$  The homogeneous coordinate ring of $\P^n$ is $S = k[x_0, \ldots, x_n].$  We let $I_X  = \oplus H^0(\P^n, \I_X(d))$ and $S_X = S/I_X.$  We let $H$ denote a general hyperplane in $\P^n$ and its pullback under a morphism.  We write $\mO(k)$ for $\mO(kH)$ when no confusion will arise.  We may write $H^i(\F)$ (resp. $h^i(\F)$) for $H^i(X, \F)$, (resp. $h^i(X, \F)$) if the meaning is clear.

Let $C$ be a smooth curve of genus $g.$  Throughout, $L$ is a very ample line bundle on $C$ embedding it as a linearly normal curve in $\P^n = \P(H^0(C, L))$ with degree $d = \deg L.$

A line bundle $L$ on a smooth curve $C$ is said to {\em separate $k$ points} if $h^0(C,L(-Z))=h^0(C,L)-k$ for all $Z\in S^kC$, where $S^kC$ is the $k$th symmetric product of $C.$  We let $\Sigma_k$ denote the variety of $(k+1)$-secant $k$-planes to $C$ and write $\Sigma$ for the variety $\Sigma_1.$

We recall the first stages of a construction of Aaron Bertram which provides the geometric framework for our results. 

\begin{thm}\cite[Theorem 1]{bertram1}\label{terracini}
Suppose $L$ separates $4$ points.  Let $g:B_1\rightarrow B_0 = \P^n$ be the blowup of $B_0$ along $C$ with $\wts$ the proper transform of $\Sigma.$  Let $h:B_2 \to B_1$ be the blowup of $B_1$ along $\wts$ and  $E_i$ be the proper transform in $B_i$ of each exceptional divisor.  We further let $f=g\of h$. 

Then $\wts \subset B_1$ is smooth and irreducible, and transverse to $E_1$, so in particular $B_2$ is smooth.  Moreover, (Terracini recursiveness) if $x\in \Sigma\setminus C$, then $f^{-1}(x) \cong \P(H^0(C,L(-2V)))$, where $V$ is the unique divisor of degree $2$ whose span contains $x$. If $x\in C,$ then $f^{-1}(x)$ is isomorphic to the blowup of $\P(H^0(C,L(-2x)))$ along the image of $C$ embedded by $L(-2x)$.
{\nopagebreak \hfill $\Box$ \par \medskip}
\end{thm}

\begin{rmk}\label{recursion}
Bertram's construction continues, blowing up up the strict transform of each $\Sigma_k$ successively, so that a fiber over a point of $C$ of the composition is $\P^{n-2}$ in which we have blown up copies of $\Sigma_i$ for $i = 0, \ldots, k-1$ and the degree of $\Sigma_0 =C$ is two less than the degree of the original embedding. We will abuse notation in the hopes of highlighting the recursive nature of the construction and denote the restriction of $E_i$ to a fiber $F$ of the composition using the notation of our setup relative to the blowing up that has occurred within $F.$  For example, if $x \in C$ and $F = (h\circ g)^{-1}(x)$, we will write $\mO_{B_2}(E_2)|_F = \mO_F(E_1),$  keeping in mind that ``$E_1 \subset F$'' is the exceptional divisor of $\P^{n-2}$ blown up at $C$ where the degree has already dropped by two.
\end{rmk}

A key point in what follows is that $\wts$ is a resolution of singularities of $\Sigma$, and is a $\P^1$-bundle over $S^2C$ in a natural way.  We summarize this relationship in Lemma \ref{S^2C tech}.  

\begin{lemma}\label{S^2C tech}
The variety $\wts \subset B_1$ is a resolution of singularities $g:\wts\rightarrow\Sigma$ with the following properties
\begin{enumerate}
\item $g_*\mO_{\wts}=\mO_{\Sigma}.$  
\item $Z:=E_1 \cap \wts \cong C \times C.$
\item The restriction $g:C\times C\rightarrow C$ is projection onto one factor.  
\item The restriction of the linear system $|2H-E_1|$ to $\wts$ yields a morphism $\pi: \wts \to S^2C$ realizing $\wts$ as a $\P^1$-bundle over $S^2C.$  The restriction of this morphism to $Z \cong C\times C$ is the canonical double cover $d:C\times C\rightarrow S^2C$. 
\item  If we define $\delta$ by $d^*\mO_{S^2C}\left(\frac{\delta}{2}\right)=\mO_{C\times C}(\Delta)$, then $d_*\mO_{C\times C}=\mO_{S^2C}\oplus\mO_{S^2C}(-\frac{\delta}{2}).$
\item If $F$ is a fiber of the $\P^1$-bundle $\pi: \wts \to S^2C$, then $\mO_F(aH-bE)=\mO_{\P^1}(a-2b)$.
\end{enumerate}
\end{lemma}

\begin{proof}
The first is \cite[3.2]{vermeireidealreg}, the second and third are \cite[3.7]{vermeireflip1}, the fourth is \cite[3.8]{vermeireflip1}.  Part (5) follows from \cite[V.22]{bpv}.  For (6), note that each fiber $F$ is the proper transform of a secant line, hence the intersection with a hyperplane is $1$, while the intersection with the exceptional divisor is $2$ (since each secant or tangent line intersects $C$ in a scheme of length two).
\end{proof}

\begin{lemma}\label{tech}
With hypotheses and notation as above:
\begin{enumerate}
\item $\Sigma \subset B_0$ is normal and is smooth away from $C$.
\item $f_*\mO_{B_2}=\mO_{B_0}$ and $R^jf_*\mO_{B_2}=0$ for $j\geq 1$.
\item $R^if_*\mO_{B_2}(-E_2)=\begin{cases}
\I_{\Sigma} & i=0\\
H^1(C,\mO_C)\tensor\mO_C & i=2\\
0 & i\neq0,2.
\end{cases}$
\item $R^ig_*\mO_{B_1}(-mE_1) = R^ih_*\mO_{B_2}(-mE_2)=0$ for $i>0$ and $m\geq 0.$
\item $R^ig_*\mI_{\wts} = R^if_*\mO_{B_2}(-E_2).$
\item $R^if_*\mO_{B_2}(-E_1-E_2)=\I_{\Sigma/\P^n}$ for $i=0$ and is zero otherwise.
\end{enumerate}
\end{lemma}

\begin{proof}
The first two can be found in \cite[3.2]{vermeireidealreg}, while the third is \cite[Proposition 9]{vermeiresecreg} and the fourth is Lemma 4.3.16 in \cite{lazarsfeld}.  Part (5) follows immediately from (4) and a degenerate case of Grothendieck's composition of functors spectral sequence \cite{tohoku}.

For the sixth item, we compute sheaves $R^if_*\mO_{E_1}(-E_2)$ and use them to show the claim via
$$\ses{\mO_{B_2}(-E_1-E_2)}{\mO_{B_2}(-E_2)}{\mO_{E_1}(-E_2)}.$$

Since $E_1 \to C$ is flat, the locally free sheaf $\mO_{E_1}(-E_2)$ is also flat over $C.$  Thus, we can compute higher direct images via cohomology along the fibers of $f$ restricted to $E_1$ by \cite[Corollary III.12.9]{hartshorne}.  By the Terracini recursiveness portion of Theorem~\ref{terracini}, if $x \in C,$ a fiber $F = f^{-1}(x)$ is the blowup of $C$ in $\P H^0(C, L(-2x))$ and $E_2$ intersects $F$ in the exceptional divisor $E_1$ of this blowup.  As $H^i(F,\mO_{F}(-E_1))=H^i(\P(H^0(C,L(-2x))),\I_C),$ it is clear that $H^i(F, \mO_F(-E_1))$ vanishes for $i=0,1,$ and $h^2(\P(H^0(C,L(-2x))),\I_C) = h^1(C,\mO_C)=g$.  We conclude that $R^if_*\mO_{E_1}(-E_2)= 0$ for $i = 0,1$ and that for $i=2$ it is locally free of rank $g.$  Note that by part (5), $R^2f_*\mO_{B_2}(-E_2)$ is also locally free of rank $g.$  Therefore, if the map between them is a surjection, it is an isomorphism.  

To get the surjectivity above we show $R^3f_*\mO_{B_2}(-E_1-E_2)=0$ by looking at
\begin{equation}\label{E_2 ses}
0\to \mO_{B_2}(-E_1-E_2) \to \mO_{B_2}(-E_1) \to \mO_{E_2}(-E_1) \to 0
.\end{equation}
Applying $h_*,$ the projection formula and the observation that $E_2 \to \wts$ is a projective bundle, we see that 
\[
0\to \mI_{\wts}(-E_1) \to \mO_{B_1}(-E_1) \to \mO_{\wts}(-E_1) \to 0
\]
is exact and all higher direct images vanish.  If we apply $g_*$ we get
\[
\to R^2g_*\mO_{\wts}(-E_1) \to R^3g_*\I_{\wts}(-E_1) \to R^3g_*\mO_{B_1}(-E_1) \to,
\]
where the first term vanishes because $\wts \to \Sigma$ has fibers of dimension at most one, and the third term vanishes by (4).
\end{proof}

We will use Lemma \ref{conormal} to show that $H^1(\Sigma, \mO_{\Sigma}(2))=0$ in Proposition \ref{gw}.
\begin{lemma}\label{conormal}
Let $L$ be a very ample line bundle on a variety $X$ with $H^i(X,L)=0$ for $i>0$, $E$ a locally free sheaf on $X$.  Let $\varphi:X\rightarrow\P^n=\P(H^0(X,L))$ be the induced morphism.  Then
\begin{enumerate}
\item $H^i(X\times X,\left(L\boxtimes E\right)\tensor\I_{\Delta})=H^i(X,\varphi^*\Omega^1_{\P^n}\tensor L\tensor E)$
\item $H^i(X\times X,\left(L\boxtimes E\right)\tensor\I^2_{\Delta})=H^i(X,N^*_{X/\P^n}\tensor L\tensor E)$
\end{enumerate}
\end{lemma}

\begin{proof}
Applying $\left(\pi_2\right)_*$ to the exact sequence
$$\ses{\left(L\boxtimes E\right)\tensor\I_{\Delta}}{L\boxtimes E}{\left(L\boxtimes E\right)\tensor\mO_{\Delta}}$$
yields a twist of the Euler sequence on $X$:
$$\ses{\varphi^*\Omega^1_{\P^n}\tensor L\tensor E}{H^0(X,L)\tensor E}{L\tensor E}$$
Note that the hypothesis $H^i(X,L)=0$ and the fact that $L$ is globally generated imply that all higher direct images vanish, and part (1) follows immediately. 

As $\mO_{\Delta}\tensor\I_{\Delta}=N^*_{\Delta}=\Omega^1_X$; applying $\left(\pi_2\right)_*$ to the exact sequence
$$\ses{\left(L\boxtimes E\right)\tensor\I^2_{\Delta}}{\left(L\boxtimes E\right)\tensor\I_{\Delta}}{\left(L\boxtimes E\right)\tensor N^*_{\Delta}}$$
yields a twist of the conormal sequence on $X$:
$$\ses{N^*_{X/\P^n}\tensor L\tensor E}{\varphi^*\Omega^1_{\P^n}\tensor L\tensor E}{\Omega_X^1\tensor L\tensor E}$$
Note that the hypothesis $H^i(X,L)=0$ and the fact that $L$ is very ample imply that all higher direct images vanish, and part (2) follows similarly.
\end{proof}

\section{$\Sigma$ is ACM}

 The main goal of this section is the proof of Theorem~\ref{intromainACM}.
As a consequence of our work we get Corollary \ref{projnor} showing that $\Sigma$ is projectively normal.  We will work throughout with the following hypothesis.

\begin{rem}{Hypothesis}{dagger}
Let $C \subset \P^n$ be a smooth linearly normal curve of genus $g$ and degree $d \geq 2g+3.$ \end{rem}

Using the Serre-Grothendieck correspondence between local and global cohomology, the depth of the maximal ideal on the homogeneous coordinate ring of $\Sigma \subset \P^n$ can be measured by vanishings of global cohomology groups.  We see that $\Sigma$ is ACM if and only if  $H^i(\P^n,\I_{\Sigma}(k))=0$ for all $k$ and for $0<i\leq\operatorname{dim} \Sigma$ (e.g. \cite[Ex. 18.16]{eisenbud}). In light of  \cite{vermeiresecreg} where it is shown that $\I_{\Sigma}$ is $5$-regular, in order to show that $\Sigma$ is ACM we are left to show that $H^i(\Sigma,\mO_{\Sigma}(k))=0$ for $i=1,2$ and all $k\leq 3-i$.  In what follows we handle the required cohomological vanishing cases individually.

\subsection{Vanishings for $k<0$}
The vanishings needed for $k<0$ follow easily from Kawamata-Viehweg vanishing together with part (3) of Lemma \ref{tech}.  We write the 5-term sequence associated to the Leray spectral sequence (applying Theorem \ref{terracini}) to the map $g: \wts \to \Sigma$ as it will be crucial in what follows (note that the first and fourth terms follow by part (1) of Lemma~\ref{S^2C tech}).

\begin{equation}\label{5-term f^1}
\begin{split}
0 \to H^1(\Sigma, \mO_{\Sigma}(k)) \to H^1(\wts, \mO_{\wts}(k)) \to H^0(\Sigma, R^1g_*\mO_{\wts}(k))\\
 \to H^2(\Sigma, \mO_{\Sigma}(k)) \to H^2(\wts, \mO_{\wts}(k))
\end{split}
\end{equation}

\begin{thm}\label{negative} If $C$ satisfies Hypothesis~\ref{dagger}, then $H^i(\Sigma,\mO_{\Sigma}(k))=0$ for $k<0$ and $i=1,2.$
\end{thm}
\begin{proof}
We know that $g^*\mO_{\Sigma}(1) =  \mO_{\wts}(1)$ is big and nef on $\wts,$ hence $H^i(\wts,\mO_{\wts}(k))=0$ for $k<0$ and $i<3$ by Kawamata-Viehweg vanishing.  Using the sequence (\ref{5-term f^1}), we have the claimed vanishing for $i=1$ immediately.  As $R^1g_*\mO_{\wts}\cong H^1(C,\mO_C)\tensor\mO_C$ by Lemma \ref{tech} (3-5), we have $H^0(\Sigma,R^1g_*\mO_{\wts}(k))=H^1(C,\mO_C)\tensor H^0(C,\mO_C(k))=0,$ and the vanishing for $i=2$ also follows.
\end{proof}

\subsection{Vanishings of $H^1(\Sigma, \mO_{\Sigma}(k))$ for $k>0$}
All of the remaining vanishings exploit the structure of $\wts$ as a $\P^1$-bundle over $S^2C.$   Given work of the second author in \cite{vermeiresecreg}, the projective normality of $\Sigma$ follows by exploiting Terracini recursion as a corollary of Proposition \ref{gw}.

\begin{prop}\label{gw}
If $C$ satisfies Hypothesis~\ref{dagger}, then $H^1(\Sigma,\mO_{\Sigma}(2))=0$.
\end{prop}

\begin{proof}
We show that $H^2(\P^n,\I_{\Sigma}(2))=0$.  

Since $\mO(2H-E)$ is trivial along the fibers of $\pi: \wts \to S^2C,$  $\mO_{\wt{\Sigma}}(2H-E)=\pi^*M$ for some line bundle $M$ on $S^2C$ \cite[Ex. III.12.4]{hartshorne}. 
From \cite[3.6]{vermeireidealreg} we know that $$\mO_{\wt{\Sigma}}(2H-E)\tensor \mO_Z\cong \pi^*M\tensor\mO_Z\cong L\boxtimes L\tensor\mO_Z(-2\Delta)$$  Further restricting $\pi$ to the double cover $d:C\times C\rightarrow S^2C$, by the projection formula and part (5) of Lemma~\ref{S^2C tech} we have $$H^i(Z,L\boxtimes L\tensor\mO_Z(-2\Delta))=H^i(S^2C,M)\oplus H^i\left(S^2C, M\tensor \mO_{S^2C}\left(-\frac{\delta}{2}\right)\right).$$  Again by the projection formula, we know that $H^i(\wts,\mO(2H-E))=H^i(S^2C,M)$.  By Lemma~\ref{conormal}, we have $H^i(Z,L\boxtimes L\tensor\mO_Z(-2\Delta))\cong H^i(C,N_C^*(2))$.  Thus we immediately have $H^2(Z,L\boxtimes L\tensor\mO_Z(-2\Delta))=0$, but this in turn implies $H^2(S^2C,M)=H^2(\wts,\mO(2H-E))=0$.

Let $\mathcal{L}_L$ be the line bundle on $S^2C$ such that $d^*\mathcal{L}_L=L\boxtimes L$ (e.g. \cite[\S 2.1]{kouv}).  Now, as $L\boxtimes L\tensor\mO_Z(-\Delta) =d^*\left(\mathcal{L}_L\tensor\mO_{S^2C}\left(-\frac{\delta}{2}\right)\right)$,  we know that
\begin{eqnarray*} 
d_*\left(\left(L\boxtimes L\right)\tensor\mO_Z(-\Delta)\right) & = & \left[\mathcal{L}_L\tensor\mO_{S^2C}\left(-\frac{\delta}{2}\right)\right]\oplus \left[\mathcal{L}_L\tensor\mO_{S^2C}\left(-2\frac{\delta}{2}\right)\right]\\
&=& \left[\mathcal{L}_L\tensor\mO_{S^2C}\left(-\frac{\delta}{2}\right)\right]\oplus M
\end{eqnarray*}
Again by Lemma~\ref{conormal} we know that $H^1(C\times C,L\boxtimes L\tensor\mO_Z(-\Delta))=H^1(C,\Omega^1_{\P^n}(2)\tensor\mO_C)=0$, where the vanishing comes from quadratic normality of the embedding of $C$.  Thus $H^1(S^2C,M)=H^1(\wts,\mO_{\wts}(2H-E))=0$.

We see immediately that $H^2(B_1,\I_{\wt{\Sigma}}(2H))=H^1(\wt{\Sigma},\mO_{\wt{\Sigma}}(2H))$, and from the sequence 
$$\ses{\mO_{\wt{\Sigma}}(2H-E)}{\mO_{\wt{\Sigma}}(2H)}{\mO_{\wt{\Sigma}}(2H)\tensor\mO_E}$$
and the (just proved) fact that $H^i(\wts,\mO_{\wts}(2H-E))=0$ for $i=1,2$ implies further that $H^2(B_1,\I_{\wt{\Sigma}}(2))=H^1(\wt{\Sigma},\mO_{\wt{\Sigma}}(2)\tensor\mO_E)$.  A straightforward computation gives:
\begin{eqnarray*}
h^1(\wt{\Sigma},\mO_{\wt{\Sigma}}(2H)\tensor\mO_E) & = & h^1(C\times C,L^2\boxtimes\mO_C)\\
 & = & h^0(C,L^2)\cdot h^1(C,\mO_C)\\
 & = & h^0(C,H^1(C,\mO_C)\tensor L^2)\\
 & = & h^0(\P^n,R^2g_*\I_{\wt{\Sigma}}(2)).\\
\end{eqnarray*}
Therefore,  $h^2(B_1,\I_{\wt{\Sigma}}(2))=h^0(\P^n,R^2g_*\I_{\wt{\Sigma}}(2)).$

Interpreting what we have just shown in terms of the Leray-Serre spectral sequence associated to $g_*\mI_{\wts}(2)$, we have $h^2(B_1, \I_{\wt{\Sigma}}(2))=\operatorname{dim}E_2^{0,2}$.  We also know that $R^1g_*\mI_{\wts}(2) = 0$ by the projection formula and Lemma \ref{tech} (3) and (5).  Thus, at the $E_2$ level, where we have
$$0\rightarrow E_2^{0,1}\stackrel{d_2}{\rightarrow}E_2^{2,0}\rightarrow 0$$
and 
$$0\rightarrow E_2^{0,2}\stackrel{d_2}{\rightarrow}E_2^{2,1}\rightarrow 0$$
we see that $E_{2}^{2,0}=E_{\infty}^{2,0}$ and $E_{2}^{0,2}=E_{\infty}^{0,2}$ because $H^i(\P^n,R^1g_*\mI_{\wts}(2))=0$.  
Recall \cite[5.2.6]{weibel} that $H^2:=H^2(B_1,\I_{\wt{\Sigma}}(2))$ has a finite filtration
$$0=F^{3}H^2\subseteq F^2H^2\subseteq F^1H^2\subseteq F^0H^2=H^2$$
where $F^2H^2\cong E_{\infty}^{2,0}$ and $H^2/F^1H^2\cong E_{\infty}^{0,2}$.  

Now, because $\operatorname{dim}H^2=\operatorname{dim}E_2^{0,2}=\operatorname{dim}E_{\infty}^{0,2}$, we have $F^1H^2=0$, but this implies that $F^2H^2=E_{\infty}^{2,0}=0$, and hence that $E_2^{2,0}=0$.
\nopagebreak \hfill $\Box$ \par \medskip
\end{proof}

In \cite{vermeiresecreg} it was shown that for the general embedding of degree at least $2g+3$, $\Sigma$ is projectively normal;  the only vanishing that could not be shown to always hold was $H^1(\P^n,\I_{\Sigma}(2))=0$.  Proposition~\ref{gw} allows us remove the hypothesis that the embedding must be general.  The idea in \cite{vermeiresecreg} was to obtain a vanishing statement for direct image sheaves, and then to use those vanishings along with \cite[p.52,Cor $1\frac{1}{2}$]{mumford} to show that the cohomology groups along the fibers vanish.  Of course, to make this work, we must find a flat morphism and a locally free sheaf so that the restriction of the sheaf to the fiber is precisely the vanishing statement we want.  This is done using Theorem~\ref{terracini}.  However, note that in the proof we need to increase the degree of the embedding to at least $2g+5$, so that curves of degree $2g+3$ occur in the fibers.

\begin{cor}\label{projnor}
Let $C\subset\P^n$ be a smooth curve embedded by a line bundle $L$ of degree at least $2g+3$.  Then $\Sigma$ is projectively normal.
\end{cor}
\begin{proof}
We know by combining \cite[Proposition 12]{vermeiresecreg} with \cite[1.16]{wahl} that $H^1(\P^n,\I_{\Sigma}(k))=0$ for $k=1,3$, and by \cite[Corollary 11]{vermeiresecreg} that $H^1(\P^n,\I_{\Sigma}(k))=0$ for $k\geq 4$.  Clearly, $H^1(\Sigma,\mO_{\Sigma}(2))=H^2(\P^n,\I_{\Sigma}(2))$.  As these vanish by Proposition~\ref{gw}, we note that by Lemma~\ref{tech} we have $H^2(B_2,\mO(2H-E_1-E_2))=0$.  We further have $H^i(B_2,\mO(2H-E_1-E_2))=0$ for $i\geq 3$ by $5$-regularity of $\I_{\Sigma}$.

Also by Lemma~\ref{tech}, along the fibers of $E_1\rightarrow C$ we are computing $H^i(\P^{n-2},\I_C(1))$, thus $R^if_*\mO_{E_1}(2H-E_1-E_2)=0$ for $i\geq0$; this gives $H^i(B_2,\mO_{E_1}(2H-E_1-E_2))=0$ and consequently that $H^i(B_2,\mO(2H-2E_1-E_2))=0$ for $i\geq 2$.

Fixing a point $p\in C$, and applying an extension of Theorem~\ref{terracini} to $L(2p)$ (which now separates $6$ points as $L$ is non-special), we may blow up three times to get a resolution of $\Sigma_2$.  In the notation of \cite[Theorem 15]{vermeiresecreg}, the previous paragraph gives $R^if_*\mO_{E_1}(kH-2E_1-2E_2-E_3)=0$ for $i\geq 2$, since the restriction of $\mO_{E_1}(kH-2E_1-2E_2-E_3)$ to a fiber of $E_1\rightarrow C$ is $\mO(2H-2E_1-E_2)$ using the convention of Remark \ref{recursion}.  It was shown in \cite[Theorem 15]{vermeiresecreg} that $R^1f_*\mO_{E_1}(kH-2E_1-2E_2-E_3)=0$, and so we know that $H^1$ along the fibers vanishes by \cite[p.52,Cor $1\frac{1}{2}$]{mumford}.  Thus we have $H^1(\wtp,\mO(2H-2E_1-E_2))=0$ and so, as above, $H^1(\wtp,\mO(2H-E_1-E_2))=H^1(\P^n,\I_{\Sigma}(2))=0$.  
\end{proof}

\begin{thm}\label{k=1} If $C$ satisfies Hypothesis~\ref{dagger}, then $H^i(\Sigma,\mO_{\Sigma}(1))=0$ for $i=1,2.$
\end{thm}
\begin{proof}
We first show that $H^i(\Sigma,\mO_{\Sigma}(1))=0$ vanishes for $i=2,$ from which the vanishing for $i=1$ follows from a computation.

Note that $h^2(\Sigma,\mO_{\Sigma}(1))=h^3(\P^n,\I_{\Sigma}(1))$ and $H^3(\P^n,\I_{\Sigma}(1))$ is isomorphic to $H^3(B_2,\mO_{B_2}(H-E_1-E_2))$ by the last part of Lemma \ref{tech}.  Using Equation (\ref{E_2 ses}) twisted by $H$, the projection formula, gives $R^ih_*(\mO_{E_2}(H-E_1)) = R^ih_*(\mO_{E_2}) \otimes \mO_{\wts}(H-E_1)$.  By part (6) of Lemma~\ref{S^2C tech} the restriction of $\mO(H-E_1)$ to the fibers of $\wts \to S^2C$ is isomorphic to $\mO_{\P^1}(-1)$, hence $h^i(\wts,\mO(H-E_1))=0$ for all $i,$ which implies that $h^i(E_2, \mO_{E_2}(H-E_1))=0.$   We therefore have $h^3(B_2,\mO_{B_2}(H-E_1-E_2)) = h^3(B_2,\mO_{B_2}(H-E_1)).$   

We see that $R^if_*(\mO_{B_2}(H-E_1))=0$ for $i \geq 1$ and $f_*(\mO_{B_2}(H-E_1))=\I_C(1)$ by \cite[1.2,1.4]{bel}.  Thus $h^3(B_2,\mO_{B_2}(H-E_1))=h^3(\P^n, \mI_C(1))=0$.

As $h^2(\Sigma,\mO_{\Sigma}(1))=0,$ the first three terms in the $5$-term sequence (\ref{5-term f^1}) with $k = 1$ form a short exact sequence.  Then $h^1(\Sigma,\mO_{\Sigma}(1))=0$  as we see the second and third terms have the same dimension:
\[
h^1(\wts,\mO_{\wts}(H))=h^1(S^2C,\mathcal{E})=h^1(C\times C,\pi_1^*\mO_{C}(1))=h^1(\mO_C) h^0(\mO_C(1))
\]
and
\[
h^0(\Sigma,R^1g_*\mO_{\wts}(H))=h^0(C,H^1(\mO_C)\tensor\mO_C(1))=h^1(\mO_C)h^0(\mO_C(1)).\]
\end{proof}

\begin{rem}{Remark}{canonical}
Note that in the case of a canonical curve, we have $$h^0(\Sigma,R^1g_*\mO_{\wts}(H))=h^1(C,\mO_C)\cdot h^0(C,\mO_C(1))=g^2$$
while
$$h^1(\wts,\mO_{\wts}(H))=h^1(\mO_C)\cdot h^0(\mO_C(1))+h^0(\mO_C)\cdot h^1(\mO_C(1))=g^2+1.$$
Therefore using the 5-term sequence~(\ref{5-term f^1}) again we see $h^1(\Sigma,\mO_{\Sigma}(1))\geq1$ (in fact, it can be shown to be equality).  Thus the secant variety to a canonical curve of Clifford index at least $3$ (e.g. the generic curve of genus $\geq7$) is \textit{never} ACM.

Note the secant variety of a canonical curve $C\subset \P^4$ is a hypersurface of degree $16$, hence is ACM, but such curves have Clifford index $\leq 2$.
\end{rem}

\subsection{Vanishings for $k=0$}
We now consider the vanishing of $H^i(\Sigma, \mO_{\Sigma})$ where $i=1,2.$  

\begin{prop}\label{hardh1}
If $C$ satisfies Hypothesis~\ref{dagger}, then $H^1(\Sigma,\mO_{\Sigma})=0.$
\end{prop}
\begin{proof}
Associated to the morphism $g:B_1\rightarrow \P^n$ we have
\begin{center}
{\begin{minipage}{1.5in}
\diagram
 & \dto & 0 \dto & 0 \dto &  \\
0 \rto& H^1(g_*\mO_{\wts}) \dto \rto & H^1(\mO_{\wts})\dto^{\alpha} \rto^{\gamma} & H^0(R^1g_*\mO_{\wts}) \dto \rto & \\
0 \rto& H^1(g_*\mO_{Z}) \rto & H^1(\mO_{Z})\rto^{\beta} & H^0(R^1g_*\mO_{Z}) \rto & 0
\enddiagram
\end{minipage}}
\end{center}
where the horizontal maps come from $5$-term exact sequences.

As $Z \cong C \times C,$ we see that the inclusion and projection in the bottom row come from the K\"unneth formula.  The map $\alpha:H^1(\mO_{\wts})\rightarrow H^1(\mO_{Z})$ is an inclusion because it is the diagonal mapping $\alpha:H^1(\mO_{S^2C})\rightarrow H^1(C,\mO_C)\oplus H^1(C,\mO_C)$ induced by the pull-back of $d:Z\rightarrow S^2C$ to $\wts.$  We conclude that the composition $\beta \circ \alpha$ is an isomorphism.  Moreover, as $H^0(R^1g_*\mO_{\wts})\rightarrow H^0(R^1g_*\mO_{Z})$ is an isomorphism, we see that $\gamma$ is an isomorphism by commutativity of the diagram.  Hence, $H^1(g_*\mO_{\wts})=H^1(\Sigma,\mO_{\Sigma})=0.$ 
 \end{proof}

\begin{prop}\label{h20}
If $C$ satisfies Hypothesis~\ref{dagger}, then $H^2(\Sigma,\mO_{\Sigma})=0.$
\end{prop}

\begin{proof}
We note that $h^i(\Sigma, \mO_{\Sigma}) = h^{i+1}(\P^n, \mI_{\Sigma})$ for $i=1,2.$  Moreover, $h^j(\P^n, \mI_{\Sigma}) = h^j(B_2, \mO_{B_2}(-E_1-E_2))$ by part (6) of Lemma \ref{tech}.  Therefore, the result follows if we can show that $h^2(B_2, \mO_{B_2}(-E_1-E_2)) = h^3(B_2, \mO_{B_2}(-E_1-E_2))$, since we know by Proposition~\ref{hardh1} that $h^2(B_2, \mO_{B_2}(-E_1-E_2))=0$.

To this end, consider the long exact sequence associated to Equation (\ref{E_2 ses}).  The result will follow if $h^2(B_2, \mO_{B_2}(-E_1)) = h^2(E_2, \mO_{E_2}(-E_1))$ is equal to $g$ and $h^1(E_2, \mO_{E_2}(-E_1)) = h^3(B_2, \mO_{B_2}(-E_1))=0$.

From the sequence $\ses{\mO_{B_2}(-E_1)}{\mO_{B_2}}{\mO_{E_1}}$ we see immediately that $h^i(B_2,\mO_{B_2}(-E_1))=g$ if $i=2$ and is zero otherwise as $R^jf_*\mO_{B_2} = 0$ for $j > 0$ from Lemma \ref{tech} (2) and $h^j(\mO_{E_1}) = h^j(\mO_C)$ for all $j.$

We compute the cohomology of $\mO_{E_2}(-E_1)$ using Equation (\ref{E_2 ses}).  Using the projection formula and part (4) of Lemma \ref{tech}, we see that $R^ih_*\mO_{E_2}(-E_1) = 0$ for $i>0.$  Thus, $H^i(\mO_{E_2}(-E_1)) \cong H^i(\mO_{\wts}(-E_1)).$  

To compute $H^i(\mO_{\wts}(-E_1)),$ observe that
\[
0 \to \pi_*\mO_{\wts}(-E_1) \to \pi_*\mO_{\wts} \to \pi_*\mO_Z \to R^1\pi_*\mO_{\wts}(-E_1) \to 0,
\]
with all remaining higher direct images vanishing by parts (2) and (4) of Lemma \ref{S^2C tech} and $\pi_*\mO_{\wts}(-E_1)=0$ by part (6). 

 As $\operatorname{Hom}_{\mO_{S^2C}}(\mO_{S^2C},\mO_{S^2C}\left( -\frac{\delta}{2} \right))$ is trivial, this gives rise to the natural inclusion
\[
 \pi_*\mO_{\wts} \cong \mO_{S^2C} \hookrightarrow \mO_{S^2C} \oplus \mO_{S^2C}\left( -\frac{\delta}{2} \right) \cong \pi_*\mO_Z,
\]
and we see that $H^i(\wts,\mO_{\wts})\hookrightarrow H^i(Z,\mO_Z).$ In fact, using the long exact sequence on $\wts,$  these inclusions imply that $H^i(Z,\mO_Z)\cong H^i(\wts,\mO_{\wts})\oplus H^{i+1}(\wts,\mO_{\wts}(-E_1))$.

As $h^1(S^2C,\mO_{S^2C})=g$ and $h^2(S^2C,\mO_{S^2C})= \binom{g}{2}$ by \cite{macdonald}, using the sequence $0\to \mO_{\wts}(-E_1) \to \mO_{\wts} \to \mO_Z \to 0$
together with the K\"unneth formula and the fact that $H^i(\wts, \mO_{\wts}) \cong H^i(S^2C, \mO_{S^2C}),$ implies that $h^2(E_2, \mO_{E_2}(-E_1))=g,$ and that $h^3(E_2, \mO_{E_2}(-E_1)) = \binom{g+1}{2}$.  Further, as $H^0(Z, \mO_Z) \cong H^0(\wts, \mO_{\wts}) \oplus H^{1}(\wts, \mO_{\wts}(-E_1))$, we see immediately that $H^{1}(\wts, \mO_{\wts}(-E_1))$ is 0.
\end{proof}

\begin{proof}[Proof of Theorem~\ref{intromainACM}]
As explained at the beginning of the section, in order to show that $\Sigma$ is ACM we are left to show that $H^i(\Sigma,\mO_{\Sigma}(k))=0$ for $i=1,2$ and all $k\leq 3-i$. 

The vanishings for $k<0$ were shown in Theorem~\ref{negative}.  The vanishing for $i=1$ and $k=0$ is Proposition~\ref{hardh1}, while $i=2$ and $k=0$ is Proposition~\ref{h20}.  Both vanishings for $k=1$ are found in Theorem~\ref{k=1}.  Finally, the vanishing for $i=1$ and $k=2$ is found in Proposition~\ref{gw}.
\end{proof}

As an immediate consequence of the proof of Proposition~\ref{h20} we get a sharpening of the regularity result of the second author in \cite{vermeiresecreg}.

\begin{cor}\label{h^4}
If $C$ satisfies Hypothesis~\ref{dagger}, then $\I_{\Sigma}$ has regularity 3 if $C$ is rational and regularity 5 otherwise.  \end{cor}

\begin{proof}
Running long exact sequence associated to Equation (\ref{E_2 ses}) in the proof of Proposition~\ref{h20} shows that $h^4(\P^n,\I_{\Sigma})=\binom{g+1}{2}$.
\end{proof}

\section{Betti Diagrams}
In this section we paint a picture of the shape of the Betti diagram of $S_{\Sigma}$ that parallels the discussion of the Betti diagram of a high degree curve in Chapter 8 of \cite{E2}. In \S \ref{sec: betti} we use the fact that $\Sigma$ is ACM to use duality and algebraic techniques to compute the extremal nontrival Betti numbers, $\beta_{1,3}$ (Corollary \ref{cubics}) and $\beta_{n-3, n+1}$ (Theorem \ref{duality}) as well as the Hilbert polynomial.  Independent of the Cohen-Macaulay property,  we prove a nonvanishing result about the length of the degree $(k+2)$ linear strand of $S_{\Sigma_k}$ using determinantal methods and Koszul homology (Proposition \ref{easy p} and Corollary \ref{beyond p}) in \S \ref{sec: strand}.

\subsection{Computing Betti numbers}\label{sec: betti}
We begin with a simple consequence of duality.  As $\Sigma$ is ACM, dualizing a resolution of $S_{\Sigma}$ and shifting by $-n-1$ gives a resolution of the canonical module, which is defined to be $\omega_{\Sigma} =\Ext^{n-3}(S_{\Sigma}, S(-n-1))= \oplus_{d \in \Z} H^0(\P^n, \omega^{\circ}_{\Sigma}\otimes L^d),$ where $\omega^{\circ}_{\Sigma}=\mathcal{E}xt^{n-3}_{\P^n}(\mO_{\Sigma},\mO_{\P^n}(-n-1))$ is the dualizing sheaf of $\Sigma.$  Therefore, the last few Betti numbers of $S_{\Sigma}$ are the first few of $\omega_{\Sigma}.$  As an immediate consequence of Corollary \ref{h^4} we see that the number of minimal generators of $\omega_{\Sigma}$ in degree 0 is $\binom{g+1}{2}$ and hence depends only on $g,$ independent of the embedding (as long as the degree is at least $2g+3$).

\begin{cor}\label{duality}
If $C$ satisfies Hypothesis~\ref{dagger}, then $\beta_{n-3, n+1}= \binom{g+1}{2}.$
\end{cor}
\begin{proof}
If $g=0,$ we know that $\beta_{n-3, n+1} = 0.$  If $g>0,$ then Corollary \ref{h^4} shows that $\reg S_{\Sigma} = 4.$  Hence, the $a$-invariant of $S_{\Sigma}$ is 0, so $h^0(\omega^{\circ}_{\Sigma}) = \beta_{0,0}(\omega_{\Sigma}) = \beta_{n-3, n+1}(S_{\Sigma}).$  By Serre duality, 
\[h^0(\P^n, \omega^{\circ}_{\Sigma}) = h^3(\P^n, \mO_{\Sigma}) =h^4(\P^n, \mI_{\Sigma})= \binom{g+1}{2}.\]
\end{proof}

Knowing $\beta_{n-3, n+1}$ allows us to compute the Hilbert polynomial of $S_{\Sigma}$ and to gather information about other Betti numbers inductively.  To begin this process, fix general linear forms $H_1, H_2, H_3, H_4 \in S.$  Let $X$ be the intersection of $\Sigma$ with the hyperplanes determined by $H_1$ and $H_2$ and $M = S_{\Sigma}/\langle H_1, H_2, H_3, H_4 \rangle.$  Using Corollary \ref{duality} we may compute the genus of $X,$ and a formula for the Hilbert polynomial of $S_{\Sigma}$ and $\beta_{1,3}$ follows.  First we gather together basic facts about $X.$

\begin{lemma}
If $C$ satisfies Hypothesis~\ref{dagger}, the variety $X$ is a smooth curve of degree $D = \binom{d-1}{2}-g$ embedded in $\P^{n-2}$ via the complete linear series associated to a line bundle $A$  and $S_X = S_{\Sigma}/\langle H_1, H_2 \rangle.$
\end{lemma}
\begin{proof}
All the statements follow immediately from the fact that $\Sigma$ is ACM.  The only thing that may not be immediate to the reader is the fact that $\operatorname{deg}(\Sigma)=\binom{d-1}{2}-g$, though this is certainly well-known to experts.  

To see this, take a generic $\Lambda=\P^{n-3}\subset\P^n$ and consider the induced projection $\pi:\P^n\dashrightarrow \P^2$.  Every point of intersection of $\Lambda$ with $\Sigma$ corresponds to a node of $\pi(C)$.  It is well-known that the number of nodes is $\binom{d-1}{2}-g$.
\end{proof}

We will denote the genus of $X$ by $G.$  To compute $G$ we compare the Hilbert function of $S_X$ to that of successive quotients by $H_1$ and $H_2.$

\begin{prop}
If $C$ satisfies Hypothesis~\ref{dagger}, the genus of $X$ is $G = \frac12(d-2)(d+2g-3).$
\end{prop}
\begin{proof}
Since $S_X$ is 4-regular,  $h^0(X, A^m) = mD-G+1$ for $m \geq 3.$  We also know that the ideal of $\Sigma$ is empty in degree less than three, since a quadric hypersurface vanishing on $\Sigma$ must vanish twice on $C$, but this is not possible since $C$ is non-degenerate.  Therefore, we can fill in the table of Hilbert functions below where each entry in the first two columns of the table is the sum of the entries directly above and to the right.
\[
\begin{array}{c|ccc}
& S_{\Sigma}/ \langle H_1, H_2 \rangle & S_{\Sigma}/ \langle H_1, H_2, H_3\rangle & M\\
\hline
2 & \binom{n}{2} & \binom{n-1}{2} & \binom{n-2}{2}\\
3 &  3D-G+1 & 3D-G+1-  \binom{n}{2} & 3D-G+1-  \binom{n}{2}-\binom{n-1}{2} \\
4  & 4D-G+1 & D & G-2D-1+\binom{n}{2}
\end{array}.\]
But computing graded Betti numbers via Koszul homology as in Proposition 2.7 in \cite{E2} shows that $\dim M_4 = \beta_{n-3, n+1} = \binom{g+1}{2}.$  Substituting $n = d-g$ and simplifying $G = 2D+1-\binom{d-g}{2}+\binom{g+1}{2}$ gives the desired result.
\end{proof}

The computation of the Hilbert polynomial $P_{\Sigma}(m)$ follows easily.

\begin{proof}[Proof of Theorem~\ref{thm: hilbert poly}]
Using Theorem 4.2 in \cite{E2}, the Hilbert polynomial and Hilbert function of $S_{\Sigma}$ agree for $m \geq \reg S_{\Sigma} + \pdim S_{\Sigma}-n \geq 4-3=1.$  Write $P_{\Sigma}(m) = \sum_{i=0}^3 \alpha_i \binom{m+i-1}{i}.$   As $X$ is gotten by cutting down by a regular sequence of two hyperplanes, $P_X(m) = P_{\Sigma}(m)-P_{\Sigma}(m-1)-P_{\Sigma}(m-2) = \alpha_3m+\alpha_2.$  Since $X$ is a curve of degree $D$ and genus, $G,$ we see that $\alpha_3 = D$ and $\alpha_2 = 1-G.$  Since the ideal of $\Sigma$ is empty in degrees 1 and 2, we see that $P_{\Sigma}(1) = n+1$ and $P_{\Sigma}(2) = \binom{n+2}{2}$ and the result follows.
\end{proof}

We compute $\beta_{1,3}$ and get a relationship on Betti numbers at the beginning of the resolution.
\begin{prop}\label{cubics}
If $C$ satisfies Hypothesis~\ref{dagger}, we have $\beta_{1,3} = \binom{n+1}{3}-(d-2)n-3g+1$ and $\beta_{2,4} = \beta_{1,4}+\beta_{1,3}(n+1)-\binom{n+4}{n}+P_{\Sigma}(4).$ 
\end{prop}
\begin{proof}
As observed above, the Hilbert polynomial and function of $S_X$ agree in degree 3 and higher.  Since $\beta_{1,3} = \binom{n+1}{3}-(S_{X})_3 = \binom{n+1}{3}-3D+G-1,$ which simplifies to the given formula.

 By Corollary 1.10 in \cite{E2} we get a formula for the Hilbert function of $S_{\Sigma}$ in terms of graded Betti numbers: 
\[
(S_{\Sigma})_m = \sum_{i \geq 0, j \in \Z} (-1)^i \beta_{i,j}\binom{n+m-j}{n}.
\]
When $m =4,$ we must have $j \leq 4$ for $\beta_{i,j}$ to contribute to the sum.  As we know that the ideal of $X$ does not contain any forms of degree $<3,$ the result follows. 
\end{proof}

\begin{rmk}
In the formula for $\beta_{2,4}$ we have an explicit formula for each term except $\beta_{1,4},$ which is the number of quartic minimal generators of $I_{\Sigma}.$  For $d\geq 2g+4$, we know $\beta_{1,4}= 0,$ as the ideal of $\Sigma$ is generated by cubics \cite{n3p}.  
\end{rmk}

Using duality, we get a similar result for the tail of the resolution.

\begin{thm}\label{thm: tail}
If $C$ satisfies Hypothesis~\ref{dagger}, the tail of the graded Betti diagram of $S_{\Sigma}$ has the form
\[
\begin{array}{c|ccc}
  & n-5  & n-4   & n-3\\
\hline 
0 & - & - & - \\
1 & - &-. & - \\
2 & *  &*  &  A\\
3 &*  & A+B+\binom{g+1}{2}\binom{n}{2}-\binom{g}{2}(n-3)(n-1)-G & C\\
4 & B & C + \binom{g}{2}(n-3)&  \binom{g+1}{2}
\end{array}.
\]
\end{thm}
\begin{proof}
Let $A = \beta_{n-3, n-1}, B = \beta_{n-5, n-1}$ and $C = \beta_{n-3, n}.$
We know that the canonical module $\omega_X$ is $\oplus_{n \in \ZZ} H^0(K_X \otimes A^n),$ where $K_X$ is the canonical line bundle of $X.$  By duality, $\beta_{i,j}(\omega_X) = \beta_{n-3-i, n-1-j}(S_{\Sigma}).$

  By Corollary 1.10 in \cite{E2} we get a formula for the Hilbert function of $\omega_X$ in terms of graded Betti numbers: 
\[\h^0(K_X\otimes A^m) = \sum_{i\geq 0, j \in \Z} (-1)^i\beta_{i,j}(\omega_X)\binom{n-2+m-j}{n-2}.\] By Serre duality and Riemann-Roch, $h^0(K_X \otimes A^{-1}) = h^1(A) =g(d-2).$
 Thus, $g(d-2) = (n-1)\binom{g+1}{2}+C-\beta_{n-4,n},$ which gives the desired satement.  The second statement follows from the equation
\[
G = \binom{g+1}{2}\binom{n}{2} -\binom{g}{2}(n-3)(n-1) +B -\beta_{n-4, n-1}+A.
\]
\end{proof}

In particular, if $g=2$, we have the following immediate corollary.
\begin{cor}\label{cor: tail}
If $C$ satisfies Hypothesis~\ref{dagger} and $g = 2,$ the tail has the form
\[
\begin{array}{c|ccc}
   & n-5  & n-4   & n-3\\
\hline 
0 & - & - & - \\
1 & - & - & - \\
2 & * &*  &  A\\
3 & * & A+B+d-5  & C\\
4 & B        & C+d-5 &  \binom{g+1}{2}
\end{array}.
\]
\end{cor}

Based on Example~\ref{intro g=2} and the following example, we expect $A=B=C=0.$

\begin{ex}\label{g=2}
Suppose $C$ is a genus 2 curve of degree 12 in $\PP^{10}.$   We use  Example (c) in \cite{EKS} to compute the ideal of the curve determinantally over the field of rational numbers in \cite{M2}.  We then used the code created to implement \cite{sidmanSullivant} to compute the least degree pieces of the ideals of the secant varietes.  Computing the degree, dimension, and projective dimension of the resulting ideals showed that we had actually  computed the secant ideals.

\begin{verbatim}
             0  1   2   3   4   5   6   7  8 9
      total: 1 43 222 558 840 798 468 147 17 2
          0: 1  -   -   -   -   -   -   -  - -
          1: - 43 222 558 840 798 468 147  8 -
          2: -  -   -   -   -   -   -   -  9 2
\end{verbatim}
 
While the Betti diagrams for $S_{\Sigma_1}$ and $S_{\Sigma_2}$ are
\begin{verbatim}
       0  1   2   3   4   5  6 7         0  1  2  3  4 5 
total: 1 70 283 483 413 155 14 3  total: 1 41 94 61 11 4
    0: 1  -   -   -   -   -  - -      0: 1  -  -  -  - -
    1: -  -   -   -   -   -  - -      1: -  -  -  -  - -
    2: - 70 283 483 413 155  - -      2: -  -  -  -  - -
    3: -  -   -   -   -   -  7 -      3: - 41 94 61  - - 
    4: -  -   -   -   -   -  7 3      4: -  -  -  -  - -
                                      5: -  -  -  -  6 -
                                      6: -  -  -  -  5 4
\end{verbatim}
\end{ex}

\subsection{The length of the first nonzero strand}\label{sec: strand}
We now turn to the consideration of a lower bound on the length of the minimal degree linear strand of the ideal of $\Sigma_k,$ essentially following Chapter 8B.2 of \cite{E2}.  In this section we will assume the following:
\begin{rem}{Hypothesis}{star} $C$ is a smooth curve of genus $g$ and degree $d$ embedded into $\PP^n$ via a line bundle $L$ that factors as $L = L_1 \otimes L_2,$ where $|L_1| = s$ and $|L_2|=t,$ with $1 \leq s \leq t.$  
\end{rem}

First note that part of the proof of Theorem 8.12 in \cite{E2} which is given in the case $k=0$ goes through for arbitrary $k$ and allows us to see easily that the  degree $k+2$ linear strand of the Betti diagram of $\Sigma_k$ has length at least $p.$
\begin{prop}\label{easy p}
Under the conditions of Hypothesis~\ref{star}, if $d \geq 2g+2k+1+p,$ then $\beta_{p, k+1+p} \neq 0.$
\end{prop}
\begin{proof}
Factor $L$ so that $\deg L_1 \geq g+k+1$ and $\deg L_2 = g+k+p.$  By Riemann-Roch $h^0(C, L_1) \geq k+2$ and $h^0(C, L_2) \geq k+p+1.$  Thus multiplication of sections gives rise to a 1-generic matrix of linear forms with at least $(k+2)$ rows and $(k+1+p)$ columns.  Delete rows and columns to get a $(k+2) \times (k+1+p)$ matrix which is still 1-generic as an equation making a generalized entry of the smaller matrix zero also makes a generalized entry of the larger matrix zero.  The maximal minors of the smaller matrix are resolved by an Eagon-Northcott complex of length $p.$  The resolution of this ideal is a subcomplex of the ideal of $\Sigma_k.$  The result follows.
\end{proof}

We can get a better lower bound by exhibiting an explicit nontrivial cycle in the Koszul homology of $S_{\Sigma_k}$ to show that $\beta_{s+t-2k-1, s+t-k}$ does not vanish.

In Theorem 8.15 in \cite{E2}, the following result is stated for $k=1$:

\begin{thm}[Theorem 8.15 in \cite{E2}]\label{thm: kcycle}
If $I \subset S$ is a homogenous ideal which contains no forms of degree less than or equal to $k,$ then $\beta_{i, i+k} \neq 0$ if and only if there exists $\gamma \in \wedge^iS^{n+1}(-i)$ of degree $i+k$ whose image under the differential of the Koszul complex is nonzero and lies in $I \otimes \wedge^{i-1}S^{n+1}(-i+1).$
\end{thm}
\begin{proof}
The proof goes through as in \cite{E2}, replacing one by $k$ everywhere.
\end{proof}

We show that Theorem 8.13 in \cite{E2} can be extended to the case of minors of arbitrary size.

\begin{thm}\label{thm: cycle}
Suppose that $A$ is an $(s+1) \times (t+1)$ matrix of linear forms with $s+1 \geq k+2.$  If the $s+t+1$ elements in the union of the entries of the zeroth row and column are linearly independent and some $(k+2)$ minor involving the zeroth row or column does not vanish, then $\beta_{s+t-2k-1, s+t-k} (S/I_{k+2})$ does not vanish. 
\end{thm}

\begin{proof}
By Theorem \ref{thm: kcycle} it suffices to construct an explicit cycle \[\gamma \in \wedge^{s+t-2k-1}S^{n+1}(-s-t+2k+1)\] of degree $s+t-k$ whose image under the differential is a nonzero element of $I_{k+1}\otimes \wedge^{s+t-2k-2}S^{n+1}(-s-t+2k+2).$  To do this we set some notation.

By our hypotheses, the matrix $A$ has the form
\[
A = \begin{pmatrix}
a_{0,0} & a_{0,1} & \cdots & a_{0,t}\\
a_{1,0} & a_{1,1} & \cdots & a_{1,t}\\
\cdot \\
a_{s,0} & a_{s,1} & \cdots & a_{s,t}
\end{pmatrix}
=
\begin{pmatrix}
x_0 & x_1 & \cdots & x_t\\
x_{1+t} & a_{1,1} & \cdots & x_{1,t}\\
\cdot\\
x_{s+t} & a_{s,1} & \cdots & a_{s,t}
\end{pmatrix}
\]
Since the $x_i$ are linearly independent they may be chosen as part of a basis for $S_1,$ and we may assume that $\partial(e_i) =x_i$ for $i = 0, \ldots, s+t.$

Let $\sigma \subset \{1, \ldots, s\}$ and $\tau \subset \{0, \ldots, t\}$ be sets of size $k+1$ and $\sigma_t$ denote the set gotten by adding $t$ to each element of $\sigma.$  Let $e_{\sigma_t, \tau}$ be the wedge product of $\{e_0, \ldots, e_{s+t}\} \backslash (\sigma_t \cup \tau)$ in the standard order.  Note that $e_{\sigma_t, \tau} \in \wedge^{s+t-2k-2}S^{n+1}.$

We define an element $\gamma$ which will serve as our nonzero cycle.  Informally, it is the signed sum of all of the $(k+1)$-minors of $A$ which do not involve the top row, each indexed by an element $e_{\sigma_t, \tau}$ in a natural way. More precisely, 
\[
\gamma = \sum_{\sigma, \tau}  (-1)^{(\sigma+\tau)+t(k+1)} \det(\sigma \mid \tau) e_{\sigma_t, \tau},
\]  
where we define $\sigma+\tau$ to be the sum of the union of the elements in $\sigma$ and $\tau$ and $\det (\sigma \mid \tau)$ is the minor of $A$ gotten by using the rows in $\sigma$ and the columns in $\tau.$

To complete the proof we need to show that the coefficients of $\partial(\gamma)$ are all of the $(k+2)$-minors of $A$ involving the zeroth row or column. The only basis elements which can have nonzero coefficients are  $e_{\sigma_t', \tau},$ where $\sigma' \subset \{1, \ldots, s\}$ and $|\sigma'| = k+2$ and $e_{\sigma_t, \tau'}$ where $\tau' \subset \{0, \ldots, t\}$ also has size $k+2.$

To understand the coefficient of $e_{\sigma_t, \tau'},$ note that there are $k+2$ basis elements $e_{\sigma_t, \tau}$ whose images under the differential could contain $e_{\sigma_t, \tau'}$ with nonzero coefficient.  Since $\partial(e_i) = x_i$ for $i = 0, \ldots, t,$ we see that the coefficient of $e_{\sigma_t, \tau'},$ is $\pm \det (\sigma_t \cup \{0\} \mid \tau')$ where the differential expands the determinant along the zeroth row.

Similarly, the coefficient of $e_{\sigma_t', \tau},$ is $\pm \det(\sigma_t' \mid\tau \cup \{0\}),$ the differential expands the determinant along the zeroth column.  (If $0 \in \tau,$ we repeat the zeroth column twice and get coefficient zero.)
\end{proof}

We have the following result analogous to Theorem 8.12 in \cite{E2}.

\begin{proof}[Proof of Theorem~\ref{beyond p}]We will construct a matrix $A$ corresponding to the factorization of $L = L_1\otimes L_2$ by choosing bases carefully as in the proof of Theorem 8.12 in \cite{E2}.  Let $B_i$ be the base locus of $L_i.$  Fix a basis $\beta_0, \ldots, \beta_t$ of $H^0(L_2)$ so that the divisor of $\beta_i$ is $B_2+D_i$ where $D_i$ and $B_2$ have disjoint support.  Let $D$ be the divisor consisting of the union of the points in the divisors determined by $\beta_0, \ldots, \beta_t$. Since $L_1(-B_1)$ is base-point free, a general element is disjoint from $D$ and from $B_1.$    Therefore we can pick a basis $\alpha_0, \ldots, \alpha_s$ so that the divisor of each $\alpha_i$ is $B_1+E_i$ where $E_i$ is disjoint from $D$ and from $B_1.$ 
 
We will show that the $s+t+1$ elements in the union of any row and any column of the corresponding matrix $A$ are linearly independent.  Without loss of generality, consider the top row and leftmost column.  We know that the elements of the column $\alpha_0\beta_0, \alpha_1\beta_0, \ldots, \alpha_s\beta_0$ are linearly independent, as are the elements $\alpha_0\beta_0, \alpha_0\beta_1, \ldots, \alpha_0, \beta_t.$  Suppose $\gamma$ is an element in the intersection of the two vector spaces with these bases.  This implies that the divisor of $\gamma$ contains the divisor of $\alpha_0$ and of $\beta_0.$  This implies that it must contain $D_0$ and $E_0$ as well as the base loci $B_1$ and $B_2.$   Since $\gamma \in H^0(L)$ and $\alpha_0\beta_0 \in H^0(L),$ then one is a scalar multiple of the other.  Therefore, we conclude that the union of the elements in the top row and first column form a set of $s+t+1$ linearly independent elements.  

As the matrix $A$ is 1-generic, we know that the ideal generated by its maximal minors has the expected codimension and hence some $(k+2)$-minor does not vanish.  Permuting rows and columns we can assume it is in the upper lefthand corner.  Since $I_{k+2} \subseteq I_{\Sigma_k},$ the result follows from Theorems \ref{thm: kcycle} and \ref{thm: cycle}.

If $\deg L \geq 2g+2k+p+1,$ then $L$ can be factored as the product of line bundles $L_1$ with degree at least $g+k+\lfloor (1+p)/2 \rfloor$ and $L_2$ with degree greater than or equal to $\deg L_1.$   If $L_1$ and $L_2$ are generic, then each has at least $k+2$ sections.

\end{proof}

\end{document}